\documentclass{amsart}
\usepackage{amsmath,amssymb,amsfonts,amsthm}
\usepackage{hyperref}

\title{Notes on cardinals that are characterizable by a complete (Scott) sentence}

\author {Ioannis  Souldatos}
\address{4001 W.McNichols, Mathematics Department, University of Detroit Mercy, Detroit, MI 48221, USA}
\email{souldaio@udmercy.edu}



\date{\today}

\newcommand{\omegaone}{\ensuremath{\omega_1}}
\newcommand{\lomegaone}{\ensuremath{\mathcal{L}_{\omega_1,\omega}}}
\newcommand{\alephs}[1]{\ensuremath{\aleph_{#1}}}
\newcommand{\alephalpha}{\alephs{\alpha}}
\newcommand{\alephomegaone}{\alephs{\omegaone}}
\newcommand{\alephalphaplus}{\alephs{\alpha+1}}
\newcommand{\beths}[1]{\ensuremath{\beth_{#1}}}

\newcommand{\bethomegaone}{\beths{\omegaone}}

\newcommand{\M}{\ensuremath{\mathcal{M}}}
\newcommand{\N}{\ensuremath{\mathcal{N}}}
\newcommand{\A}{\ensuremath{\mathcal{A}}}
\newcommand{\B}{\ensuremath{\mathcal{B}}}

\newcommand{\C}{\ensuremath{\mathcal{C}}}

\newcommand{\mnneat}{(\M,\N)-neat\;}
\newcommand{\mnhappy}{(\M,\N)-happy\;}

\newcommand{\mnfull}{(\M,\N)-full\;}

\newcommand{\lang}[1]{\ensuremath{\mathcal{L}_{#1}}}

\newcommand{\ch}{\ensuremath{\mathcal{CH}_{\omega_1,\omega}}}
\newcommand{\homch}{\ensuremath{\mathcal{HCH}_{\omega_1,\omega}}}
\newcommand{\ltok}{\ensuremath{\lambda^{\kappa}}}
\newcommand{\ltoomega}{\ensuremath{\lambda^{\omega}}}

\newtheorem{theorem}{Theorem}[section]
\newtheorem{lemma}[theorem]{Lemma}
\newtheorem{corollary}[theorem]{Corollary}
\theoremstyle{definition}
\newtheorem{definition}[theorem]{Definition}
\newtheorem{claim}{Claim}
\newtheorem{subclaim}{Subclaim}

\theoremstyle{observation}

\newtheorem{conjecture}{Conjecture}

\theoremstyle{remark}

\numberwithin{equation}{section}

\begin{document}

\begin{abstract} This is part I of a study on cardinals that are characterizable by Scott sentences. Building on \cite{HjorthsKnightPaper}, \cite{MalitzsHanfNumber} and
\cite{BaumgartnersHanfNumber} we study which cardinals are characterizable by a
Scott sentence $\phi$ , in the sense that $\phi$ characterizes $\kappa$,
if  $\phi$ has a model of size $\kappa$, but no models of size $\kappa^+$.

We show that the set of cardinals that are characterized by a Scott sentence is closed under successors, countable unions and countable products (cf. theorems \ref{cluster},  \ref{ltoomegahom}, and corollary \ref{chhomchcountableproducts}). We also prove that if $\alephalpha$ is characterized by a Scott sentence, at least one of $\alephalpha$ and $\alephalphaplus$ is \emph{homogeneously characterizable} (cf. definition \ref{defhch}  and theorem \ref{twohomcases}). Based on  Shelah\rq{}s \cite{ShelahBorelSquares}, we give  counterexamples that characterizable cardinals are not closed under predecessors, or cofinalities.
\end{abstract}

\maketitle

\textbf{Acknowledgment:} I would like to thank the Department of Mathematics and Statistics of the University of Melbourne, Australia, for their kind hospitality for the whole academic year 2006-2007. This paper was written while visiting  professor Greg Hjorth\footnote{Professor Hjorth died unexpedently in Janury 2011.}, my thesis advisor at the above university.

\section{Introduction and known results}
This section contains the basic definitions and background theorems. A similar discussion also appears in \cite{LinearOrderings}.

Let the signature of our logic be \lang{}. We will consider only countable \lang{}. For basic definitions in infinitary logic \lomegaone,  the reader can refer to \cite{KeislersModelTheory}. In $\lomegaone$ we allow formulas that have negation, universal/existential quantification and countably long disjunctions and conjunctions, but not countably long quantification.

\begin{definition} A sentence $\sigma\in \lomegaone$ is called \emph{complete}, if for every sentence $\tau\in \lomegaone$,
either $\sigma\Rightarrow\tau$ is valid, or $\sigma\Rightarrow\neg\tau$ is valid.
\end{definition}

Obviously, all complete sentences that hold inside the same model
are equivalent to one another. D. Scott (in \cite{ScottsDenumerablyLongFormulas}) proved
that for every countable model $\mathcal{M}$, there is a sentence
$\phi_\mathcal{M}\in\lomegaone$, such that if $\mathcal{N}$ is a countable model that also
satisfies $\phi_\mathcal{M}$, then $\mathcal{M}$ and $\mathcal{N}$
are isomorphic. $\phi_\mathcal{M}$ is called the \emph{Scott
sentence} for $\mathcal{M}$.  From this, it follows that every
Scott sentence is a complete sentence.

Compactness fails in $\lomegaone$ and we also lose the upward
Lowenheim-Skolem theorem. So, it can be the case that a certain
sentence has models in cardinality \alephalpha, but not in
cardinality \alephalphaplus. This motivates the following definition:

\begin{definition}\label{charcarddef} We say that a \lomegaone-sentence $\phi$ \emph{characterizes}
\alephalpha, or that $\alephalpha$ is\emph{ characterizable}, if $\phi$ has models in all cardinalities up to $\alephalpha$, but not in cardinality
\alephalphaplus. If $\phi$ is the Scott sentence of a
countable model (or any other complete sentence), we say that it
\emph{completely characterizes} \alephalpha, or that \alephalpha\; is
\emph{completely characterizable}. Moreover, if $\phi$ is the Scott sentence of a
countable model \M, we also say that \M\; characterizes \alephalpha.

Denote by \ch, the set of all
completely characterizable cardinals.
\end{definition}
Note that the downward Lowenheim-Skolem theorem still holds, which
means that every sentence that characterizes some cardinal $\alephalpha$ has models in all cardinalities $\le\alephalpha$. 

On the other hand, W. Hanf was the first one to notice that there exists a cardinal,
call it  $\mathcal{H}(\lomegaone)$, such that, if an $\lomegaone$- sentence has a
model of this cardinality, then it has models of \emph{all}
cardinalities (see \cite{HanfsIncompactness}). $\mathcal{H}(\lomegaone)$ is called
the Hanf number for \lomegaone and it is proven to be equal to
\bethomegaone. So, $\ch\subset\bethomegaone$ and from now on we only consider cardinals that live below
\bethomegaone. We will also restrict ourselves to cardinals that are completely characterizable, and we may refer to them as just characterizable cardinals.

Using a stronger notion of characterizability which we will call \emph{homogenous characterizability}, Malitz (cf.\cite{MalitzsHanfNumber}) and Baumgartner (cf. \cite{BaumgartnersHanfNumber}) proved that for all $\alpha<\omega_1$, $\beths{\alpha+1}$ is homogeneously characterizable. We give the definition first:
\begin{definition}\label{defhch} If $P$ is a unary predicate symbol, we say that it is
\emph{completely homogeneous} for the \lang{}- structure \A, if
$P^{\A}=\{a|\A\models P(a)\}$ is infinite and every permutation of
it extends to an automorphism of \A.

If $\kappa$ is a cardinal, we will say that $\kappa$ is
\emph{homogeneously characterizable} by
$(\phi_{\kappa},P_{\kappa})$, if $\phi_{\kappa}$ is a complete
\lomegaone- sentence and  $P_{\kappa}$ a unary predicate in the
language of $\phi_{\kappa}$ such that
\begin{itemize}
    \item $\phi_{\kappa}$ does not have models of power $>\kappa$,
    \item if \M\; is the (unique) countable model of
    $\phi_{\kappa}$, then $P_{\kappa}$ is infinite and completely
    homogeneous for \M\; and
    \item there is a model \A\; of $\phi_{\kappa}$ such that
    $P_{\kappa}^{\A}$ has cardinality $\kappa$.
\end{itemize}

 If $(\phi_{\kappa},P_{\kappa})$ characterize $\kappa$ homogeneously, write
 $(\M,P(\M))\models (\phi_{\kappa},P_{\kappa})$ for that. Denote the  set of all homogeneously
characterizable cardinals by \homch. 
\end{definition}
Obviously, $\homch\subset\ch$, but the inverse inclusion fails, with $\aleph_0$ being a counterexample (cf. \cite{HjorthsKnightPaper}). It is open whether there is any other counterexample or not. By corollary  \ref{kappatoomega} it is consistent that all such counterexamples must have cofinality $\omega$. Our conjecture is that a characterizable cardinal is not homogeneously characterizable iff it
has cofinality $\omega$ (cf. conjecture \ref{cofomega}).

In \cite{MalitzsHanfNumber}, J. Malitz   proved   that under the assumption of
GCH, for every successor $\alpha<\omegaone$, $\beths{\alpha}$ is
homogeneously characterizable. J. Baumgartner improved this result
in \cite{BaumgartnersHanfNumber}  by eliminating the GCH assumption.

G. Hjorth in \cite{HjorthsKnightPaper} extended a result of Julia Knight (cf. \cite{KnightsCompleteSentence}) that
$\alephs{1}$ is characterizable, to all $\alephs{\alpha}$'s being
characterizable, for  $\alpha$ countable. 

Breaking the arguments down we see that we can easily get the following generalizations: Baumgartner's argument proves that the class of homogeneously
characterizable cardinals is closed under the powerset operator, i.e.
if $\alephs{\alpha}\in\homch$, then $2^{\alephs{\alpha}}\in\homch$ (cf. theorem \ref{baumgartner}).
Hjorth\rq{}s argument proves that the class of characterizable cardinals is closed under
successors and countable unions, i.e. if $\alephalpha\in\ch$ and
$\beta<\omegaone$, then $\alephs{\alpha+\beta}\in\ch$ (cf. theorem \ref{cluster}). This means
that characterizable cardinals come into clusters of length
$\omegaone$.

\begin{definition}\label{headofcluster} A cardinal $\alephs{\alpha}\in\ch$ is called the \emph{head of a
cluster}, if we can not find ordinals $\beta,\gamma$ such that
\begin{itemize}
    \item $\alephs{\gamma}\in\ch$,
    \item $\beta<\omegaone$ and
    \item $\alephalpha=\alephs{\gamma+\beta}$
\end{itemize}
\end{definition}

It is immediate that all characterizable cardinals are of the form
$\alephs{\alpha+\beta}$, where $\alephs{\alpha}$ is the head of a
cluster and $\beta<\omegaone$.

Since not all characterizable cardinals are homogeneously characterizable, the theorems by Hjorth and Baumgartner can not be combined directly.

Our contributions: Malitz's proof that $2^{\alephs{0}}$ is (homogeneously)
characterizable generalizes to the following: if $\lambda$ is characterizable,
then $\ltoomega$ is homogeneusly characterizable. Using this theorem, we prove  closure under countable products for both the  characterizable and homogeneously characterizable cardinals (corollary \ref{chhomchcountableproducts}). Moreover, if $\alephalpha^{\alephs{\beta}}\in\ch$, then for all
$\gamma<\omegaone$, $\alephs{\alpha+\gamma}^{\alephs{\beta}}\in\homch$ (theorem \ref{headpower}). This
means that the powers of the head of a cluster determine the
behavior of the whole cluster. We then conclude that if \C\; is the
smallest set of characterizable cardinals that contains $\alephs{0}$
and is closed under successors, countable unions, countable products and powersets, then
it is also closed under powers. This is theorem \ref{Cclass}.
Whence, we see that the class of characterizable and homogeneously characterizable cardinals
 maybe much richer that just containing the countable aleph and beth numbers. Of course, it depends on our set-theoretic universe. 

In the last section, we provide counterexamples that characterizable cardinals are not closed under predecessors and cofinalities.

\textbf{Structure of the Paper}
\begin{itemize}
  \item In section \ref{successors}, most of the theorems either follow, or extend, theorems from \cite{HjorthsKnightPaper}. The main theorem is that for every $\kappa\in\ch$, at least one of $\kappa$ or $\kappa^+$ is in $\homch$. We will use some of the theorems from this section in Part II too.
  \item In section \ref{powers} the main theorem is that for $\kappa\in\ch$, $\kappa^\omega$ is in $\homch$. The construction behind it is given in theorem \ref{ltoomegahom} and is very similar to the construction found in \cite{MalitzsHanfNumber}. Some consequences of this theorem are proved too.
  \item The last section contains counterexamples that characterizable cardinals are not closed under predecessors and cofinalities.
\end{itemize}

\section{Successors}\label{successors}

In this section we deal with successors of characterizable cardinals. The main theorem is theorem \ref{twohomcases}: If $\kappa\in\ch$, then at least one of $\kappa,\kappa^+$ is in $\homch$. 

The first two theorems are (essentially) in \cite{HjorthsKnightPaper}, although Hjorth is interested only in the case where $\alpha<\omega_1$.
\begin{theorem}[Hjorth] \label{successorthrm} If $\alephs{\alpha}\in\ch$, then $\alephalphaplus\in\ch$.
\begin{proof} Follows by (the proof of) theorem 5.1 from \cite{HjorthsKnightPaper}.
\end{proof}
\end{theorem}

\begin{theorem}[Hjorth]\label{limitthrm} Whenever $\alephs{\alpha_n}$, $n\in\omega$, is an non-decreasing sequence of cardinals in \ch,  then $\alephs{\lambda}= \sup{\alephs{\alpha_n}}$ is also in \ch.
\begin{proof} Take the disjoint union of structures that characterize $\alephs{\alpha_n}$, for all $n$.
\end{proof}
\end{theorem}

Combining these two theorems we get by induction on $\beta$:
\begin{theorem}\label{cluster} If $\alephs{\alpha}\in\ch$, then $\alephs{\alpha+\beta}\in\ch$, for $\beta<\omegaone$.
\end{theorem}
So, characterizable cardinals come into clusters of length \omegaone.

Next, for the shake of completeness, we repeat some definitions from \cite{HjorthsKnightPaper}, as well as a corollary and a lemma.

\begin{definition}Let $\M,\N$ be structures with languages $\lang{\M},\lang{\N}$ respectively. Assume that $\lang{\M}$ and $\lang{\N}$ have no common symbols and are entirely relational. Let $S$ be a ternary relation, $P$ a binary relation and for every $k\in\omega$, $T_k$ be a $(k+2)$-ary relation, all of which are new symbols and do not appear in $\lang{\M}$ or $\lang{\N}$. Let $\lang{}(\lang{\M},\lang{\N})$ be the language generated by $\lang{\M},\lang{\N}, S, P, \{T_k:k\in\omega\}$. Define:

A $\lang{}(\lang{\M},\lang{\N})$ structure $\A$ including $\M,\N$ is \mnneat if
\begin{description}
  \item[(1)] $\A\setminus (\M\cup\N)$ is finite;
  \item[(2)] for any relation symbol $R$ in $\lang{\M}$ and $\textbf{a}\in \A$, if $\A\models R(\textbf{a})$, then $\textbf{a}\in\M$;
  \item[(3)] similarly for any relation symbol $R$ in $\lang{\N}$ and $\textbf{a}\in \A$, if $\A\models R(\textbf{a})$, then $\textbf{a}\in\N$;
  \item[(4)] if $\A\models T_k(a_0,a_1,\textbf{b})$, then $a_0,a_1,\textbf{b}\in\A\setminus (\M\cup\N)$;
  \item[(5)] if $\A\models P(a,c)$, then $a\in\N$ and $c\notin (\M\cup\N)$;
  \item[(6)] if $\A\models S(a_0,c_0,c_1)$, then $a\in\M$ and $c_0,c_1\notin (\M\cup\N)$;
  \begin{description}
   \item[($*1^{\M,\N}$)] $\A$ satisfies the conjunction of:
                \begin{description}
                  \item[($*1_a^{\M,\N}$)] $\forall c_0,c_1 \bigvee_{a\in\M} S(a,c_0,c_1)$
                  \item[($*1_b^{\M,\N}$)] $\forall c_0,c_1 \bigwedge_{a_0\neq a_1} S(a,c_0,c_1)\Rightarrow \neg S(a_1,c_0,c_1)$;
                \end{description}
    \item[($*2^{\M,\N}$)] $\A$ satisfies the conjunction of:
                \begin{description}
                  \item[($*2_a^{\M,\N}$)] for each $k$ and $i<j<k$

                  $\forall c_0,c_1,b_0,\ldots,b_{k-1}\; T_k(c_0,c_1,b_0,\ldots,b_{k-1})\Rightarrow$

                   $(b_i\neq b_j \wedge c_0\neq c_1)$
                  \item[($*2_b^{\M,\N}$)] for each $k$ and permutation $\pi$ of k

                  $\forall c_0,c_1,b_0,\ldots,b_{k-1}\; T_k(c_0,c_1,b_0,\ldots,b_{k-1})\Leftrightarrow$

                  $T_k(c_1,c_0,b_{\pi(0)},\ldots,b_{\pi(k-1)})$
                  \item[($*2_c^{\M,\N}$)] $\forall c_0,c_1,b_0,\ldots,b_{k-1}\;\bigwedge_{a\in\M} \wedge_{i<k}\; T_k(c_0,c_1,b_0,\ldots,b_{k-1})\Rightarrow (S(a,c_0,b_i)\Leftrightarrow S(a,c_1,b_i))$
                  \item[($*2_d^{\M,\N}$)] $\forall c_0,c_1,b_0,\ldots,b_{k-1}\forall d \;\bigwedge_{a\in\M} (T_k(c_0,c_1,b_0,\ldots,b_{k-1})\wedge_{i<k} d\neq b_i)\Rightarrow$ $(S(a,c_0,d)\Rightarrow \neg S(a,c_1,d))$
                \end{description}
    \item[($*3^{\M,\N}$)] $\A$ satisfies the conjunction of:
                 \begin{description}
                  \item[($*3_a^{\M,\N}$)] $\forall c\bigvee_{a\in\N} P(a,c)$
                  \item[($*3_b^{\M,\N}$)] $\forall c \bigwedge_{a_0\neq a_1} \neg(P(a_0,c)\wedge P(a_1,c))$;
                \end{description}
    \end{description}
  \item[(7)] for all $c_0,c_1\in \A\setminus (\M\cup\N)$ there is $k\in\omega$ and $\textbf{b}\in\A$ with
  \[\A\models T_k(c_0,c_1,\textbf{b}).\]
\end{description}

\end{definition}

\begin{definition}
  An $\lang{}(\lang{\M},\lang{\N})$ structure $\A$ including $\M,\N$ is \mnhappy if for all finite $A\subset\A$ there is some \mnneat substructure $\B$ of $\A$ with $\B\supset A$.
\end{definition}

\begin{definition}  An $\lang{}(\lang{\M},\lang{\N})$ structure $\A$ is \mnfull if
        \begin{enumerate}
          \item it is \mnhappy;
          \item for all \mnneat $F\subset\A$ and \mnneat $H_0\supset F$ there is \mnneat $H_1\subset\A$ with $F\subset H_1$ and there is $i:H_0\cong H_1$ with $i|_F$ the identity.
        \end{enumerate}
\end{definition}

The following two lemmas are Corollary 5.2 and Lemma 5.4 from \cite{HjorthsKnightPaper}. 

\begin{lemma}[Hjorth]\label{mnfullinkappa}Let $\gamma\le\kappa$ be cardinals, $\kappa$ infinite. If $\M$ is a structure of size $\kappa$ and $\N$ a structure of size $\gamma$, then there exists an \mnfull structure of cardinality $\kappa$.
\end{lemma}

\begin{lemma}[Hjorth]\label{nomnfullinkappaplusplus} Under the same assumptions for $\M$ and $\N$, there is no \mnfull structure of size greater than $\kappa^+$.
\end{lemma}

\begin{theorem} \label{twohomcases} If $\kappa\in\ch$, then one of the following is the case:
\begin{enumerate}
  \item $\kappa^+\in\homch$ or,
  \item $\kappa\in\homch$.
  \end{enumerate}
\begin{proof} The $\phi$ be a complete sentence that witnesses $\kappa\in\ch$, \M\; be a model of $\phi$ of size $\kappa$ and $\N$ be a structure of size $\kappa$. By lemma \ref{mnfullinkappa}, there is a \mnfull structure of size $\kappa$ and by lemma \ref{nomnfullinkappaplusplus}, there is no \mnfull structure of size greater than $\kappa^+$. So, the proof splits into two cases\footnote{It seems to me that for a particular $\kappa$, it is independent of ZFC whether Case I or Case II holds, but I do not have a proof for that.}:
\begin{description}
  \item[Case I] There is no $(\M,\N)$- full structure of size $\kappa^+$, in which case we get that $\kappa\in\homch$.
  \item[Case II] There is a $(\M,\N)$- full structure of size $\kappa^+$, in which case there is also a $(\M,\N)$- full structure of size $\kappa^+$ where $\N$ has size $\kappa^+$. This gives $\kappa^+\in\homch$.
\end{description}
\end{proof}
\end{theorem}

Note here that although Case I and Case II in the proof of the previous theorem are exclusive the one to the other, Cases 1 and 2 in the statement of the previous theorem need not be exclusive the one to the other.

\begin{lemma}\label{nok+structure} If $\M,\N,\kappa$ are as above and $\kappa^{\omega}=\kappa$, then there is no \mnfull structure of size $\kappa^+$. \footnote{The idea of this proof was communicated to the author by professor Magidor.}
\begin{proof} Assume there is such \A\; that is \mnfull and has size $\kappa^+$. Then for every $a\in\A\setminus\M$, let $f_a$ be the function given by \[f_a(b)=m\mbox{ iff } S(m,a,b).\] This defines a family of function $\mathcal{F}=\{f_a|a\in\A\setminus\M\}$ that has size $\kappa^+$. Now, let $\A_0\subset\A\setminus\M$ be subset of cardinality $\omega$ and consider the restrictions $f_a|_{\A_0}$. By assumption there are $\kappa^\omega=\kappa$ many possible distinct such functions. Therefore, there have to be $\kappa^+$ many functions from $\mathcal{F}$ that they agree on $\A_0$. But by ($*2^{\M,\N}$) any two $f_a\neq f_b$ have to agree only on a finite set. Contradiction.
\end{proof}
\end{lemma}

\begin{theorem}\label{kappatoomega} If $\kappa\in\ch$ and $\kappa^{\omega}=\kappa$, then $\kappa\in\homch$.
    \begin{proof} By the previous lemma and Case I of (the proof of) theorem \ref{twohomcases}.
\end{proof}
\end{theorem}

It is consistent that all limit cardinals of uncountable cofinality satisfy   $\kappa^{\omega}=\kappa$ (under GCH for instance). By the above theorem, if $\kappa\in\ch$, then $\kappa\in\homch$. Hence, it is consistent that the only case that a limit cardinal $\kappa$ is in $\ch\setminus\homch$ is when $\kappa$ has cofinality $\omega$.
\begin{conjecture}\label{cofomega} For an infinite cardinal $\kappa$, $\kappa\in\ch\setminus\homch$ iff $cf(\kappa)=\omega$.
\end{conjecture}

\section{Powers}\label{powers}

Here we investigate powers of the form $\lambda^{\kappa}$, where
$\lambda,\kappa$ are characterizable. The main theorem is theorem \ref{ltoomegahom}: If $\lambda$ is in $\ch$, then $\ltoomega$ is in $\homch$. The idea behind the construction is similar to Malitz's proof that $2^{\omega}\in \homch$ in \cite{MalitzsHanfNumber}. 

\begin{theorem}[Baumgartner - \cite{BaumgartnersHanfNumber}]\label{baumgartner}
If $\kappa\in\homch$, then $2^{\kappa}\in\homch$.
\end{theorem}
Obviously, if $\lambda\le\kappa$, then, also, $\lambda^{\kappa}\in\homch$.

It is also easy to prove by induction that
\begin{lemma}
For every $\alpha<\omegaone$, there is a sentence $\sigma_{\alpha}$ such that
\[(M;<)\models\sigma_{\alpha}\mbox{ iff } (M;<)\cong (\alpha;\in). \]
\end{lemma}

The first goal is to prove
\begin{theorem}\label{ltoomegainch} If $\lambda\in\homch$, then $\ltoomega\in\ch$.
\begin{proof}

First notice that $\lambda$ is homogeneously characterizable, while $\ltoomega$ is just characterizable. Since $\lambda\in\homch$, there exists a complete sentence $\phi_{\lambda}$ in a language that contains a predicate symbol $P_{\lambda}$ such that $(\phi_{\lambda},P_{\lambda})$ homogeneously characterize $\lambda$ (cf. definition \ref{defhch}). I.e.
\begin{itemize}
    \item $\phi_{\lambda}$ does not have models of power $>\lambda$,
    \item if \M\; is the (unique) countable model of
    $\phi_{\lambda}$, then $P_{\lambda}$ is infinite and completely
    homogeneous for \M\; and
    \item there is a model \A\; of $\phi_{\lambda}$ such that
    $P_{\lambda}^{\A}$ has cardinality $\lambda$.
\end{itemize}

Now let $\lang{}$ be a signature that contains the symbols $K(\cdot)$, $F(\cdot)$, $V(\cdot,\cdot)$, $R(\cdot,\cdot)$, $E(\cdot,\cdot,\cdot)$ and $<$. The idea is to build a rooted tree of height $\omega$ where at every level we allow at most $\lambda$- splitting. $V(\cdot,\cdot)$ captures the set of vertices, with $V(n,\cdot)$ being the set of vertices of height $n<\omega$. For vertices $x,y$, $R(x,y)$ holds iff $x$ and $y$ are adjacent vertices, with $y$ being a descendant of $x$ (i.e. there is some $n$ such that $V(n,x)$, $V(n+1,y)$ and $x,y$ are connected). $F(\cdot)$ will capture the set of maximal branches through the tree. If $F(f)$, we will think $f$ as a function with domain $\omega$ and $f(n)$ will be in $V(n, \cdot)$. $E(f,n,y)$ indicates that $f(n)=y$, and we will just write $f(n)=y$ for short. 

The difficulty is to express $\lambda$- splitting, and this is where we the use the full power of the fact that $\lambda$ is in $\homch$. To every vertex $y$,  we assign a structure $M_y=M(y,\cdot)$ and we stipulate that $M_y$ together with the set of all the descedants of $y$ satisfy $\phi_\lambda$. Thus, they must have size $\le\lambda$. It takes some argument to prove that this yield a complete sentence (claim \ref{isoclaim}). 

Consider the conjunction of the following sentences:
\begin{enumerate}
  \item $(K;<)\cong(\omega;\in)$. This we can say by using the previous lemma.

  So, we will freely write $0,1,\ldots,n,n+1,\ldots$ for the elements of K.
  \item Let $V=\bigcup_{n\in K} V(n,\cdot)$. Then  $K\cup F \cup V \cup_{y\in V} M(y,\cdot)$ partition the whole space. All of them are infinite, except $V(0,\cdot)$.

  \item
    \begin{enumerate}
        \item $V(n,y)$ implies that $n\in K$. Write $y\in V(n)$ for $V(n,y)$.
        \item For $n\neq m\in K$, $V(n)\cap V(m)=\emptyset$.
        \item $V(0)=\{a\}$, where $a$ can be anything. This will be the root of the tree.
     \end{enumerate}
\item
    \begin{enumerate}
     \item If $F(f)$, we will write $f\in F$. 
    \item If $E(f,n,y)$, then $f\in F$, $n\in K$ and $y\in V(n)$. We will write $f(n)=y$ instead of $E(f,n,y)$.
\item $\forall n\in K \forall f\in F\exists! y\in V(n), f(n)=y$.  I.e. every $f\in F$ is a function with domain $\omega$ and such that $f(n)\in V(n)$, for every $n$.
    \end{enumerate}
  \item
    \begin{enumerate}
        \item $\forall y,z  (R(y,z)\Rightarrow \exists! n\in K (y\in V(n)\wedge z\in V(n+1)))$.
        \item For all $y$, the set $\{z|R(y,z)\}$ is infinite.
        \item $\forall n \forall z\in V(n+1) \exists! y R(y,z)$. By (a), this $y$ must be in $V(n)$.
        \item $\forall y_1\neq y_2\forall z (R(y_1,z)\Rightarrow \neg R(y_2,z))$, i.e. $R(y_1,\cdot)$ and $R(y_2,\cdot)$ are disjoint.
     \item $\forall n\in K\forall  f\in F (R(f(n),f(n+1)))$.
  
        \item $\forall n\in K\forall y_0,y_1,\ldots,y_n, (y_0=a \wedge_{i\le n} R(y_i,y_{i+1}))\Rightarrow \exists^{\infty} f\in F (\wedge_{i\le n} f(i)=y_i))$, i.e. every finite ``branch" can be extended to a maximal branch in infinitely many ways.
    \end{enumerate}
    \item
    \begin{enumerate}
        \item $M(y,x)$ implies that $y\in \cup_{n} V(n)$, while $x\notin \cup_{n} V(n)\cup F\cup K$.
        \item For every $y_1,y_2\in \cup_{n} V(n)$, $M(y_1,\cdot)\cap M(y_2,\cdot)=\emptyset$, i.e. the structure $M_y=M(y,\cdot)$ associated with every $y$ is unique.
        \item $\forall y$ $M(y,\cdot)\cup R(y,\cdot)\models \phi_\lambda$
        \item $\forall y\forall x\in M(y,\cdot)\cup R(y,\cdot)$, $(P_\lambda(x)\Leftrightarrow R(y,x))$.

        The above two sentences express the fact that $M(y,\cdot)$ together with $R(y,\cdot)$ satisfy $\phi_\lambda$ and $R(y,\cdot)$ is the part of the model given by the homogeneous predicate $P_\lambda$. This restricts the size of $R(y,\cdot)$ to at most $\lambda$. Since $R(y,\cdot)$ is the homogeneous part of the model, in the countable case, every permutation of it can be extended to an automorphism of the whole model. We will use this in what follows.
    \end{enumerate}

\end{enumerate}

The goal now is to show  that a structure that satisfies $(1)-(6)$ characterizes $\ltoomega$.

\begin{claim}\label{ltoomegabound} If $\M\models (1)-(6)$, then $|\M|\le \ltoomega$.
\begin{proof} First we prove by induction on $n\in K$  that $|V(n)|\le \lambda^n$.

For $n=0$, $|V(0)|=|\{a\}|=1=\lambda^0$. Assume that $|V(n)|\le \lambda^n$ and let $z\in V(n+1)$. By 5(a), there is a unique $y\in V(n)$ such that $R(y,z)$. So, \[V(n+1)=\bigcup_{y\in V(n)} R(y,\cdot).\] By 5(d), all these $R(y,\cdot)$ are disjoint and by (6) all of them have size $\le\lambda$. Thus,
\[|V(n+1)|=\sum_{y\in V(n)} |R(y,\cdot)|\le \sum_{y\in V(n)}\lambda=|V(n)|\cdot\lambda=\lambda^{n+1},\] as we want. Therefore $|\bigcup_{n\in K} V(n)|=\lambda^{<\omega}=\lambda$.

Since, for every $f\in F$, $f$ is a function from $K$ to $\bigcup_{n\in K} V(n)$, $|F|\le\ltoomega$.
Also, for every $y$, $M(y,\cdot)$ has size at most $\lambda$. So, put all together,
\[|\M|\le |K|+|F|+|\bigcup_{n\in K} V(n)|+|\bigcup_{y} M(y,\cdot)|\le\omega+\ltoomega+\lambda+\lambda\cdot\lambda=\ltoomega.\]
\end{proof}
\end{claim}

\begin{claim} There is $\M\models (1)-(6)$ and $|\M|=\ltoomega$.
\begin{proof} Take the full $\lambda$-tree $(T,R)$ of height $\omega$, with $f\in F$ being its maximal branches. The rest follows.
\end{proof}
\end{claim}

\begin{claim} \label{isoclaim}If $\M_1,\M_2$ are both countable models of (1)-(6), then there is an isomorphism $i:\M_1\cong\M_2$.
\begin{proof} First of all we observe that $(K(\M_1);<)\cong(\omega;\in)\cong(K(\M_2);<)$, so that we do not have to worry about $K$.

\begin{subclaim}\label{3} For all $n\in K$ and for all $y\in V(n)$, there is $f\in F$ with $f(n)=y$.
\begin{proof}
By 5(c) and by induction on $n$, if $y\in V(n)$, there is a (unique) sequence $y_0=a,y_1,\ldots, y_{n-1},y_n=y$ such that $\bigwedge_{i< n} R(y_i,y_{i+1})$. By 5(f) there is some $f$ that ``extends" this sequence, i.e. for all $i\le n$, $f(i)=y_i$. In particular, $f(n)=y$.
\end{proof}
\end{subclaim}

\begin{subclaim}For all $f_1\neq f_2\in F$, there exists $n>0$ such that $\forall m<n (f_1(m)=f_2(m))$ and $\forall m'\ge n (f_1(m')\neq f_2(m'))$.  We call this n the \emph{splitting point} of $f_1,f_2$, $n=s.p(f_1,f_2)$.
\begin{proof}
Since $f_1\neq f_2$, there exists $n'$ so that $f_1(n')\neq f_2(n')$.  Let $n$ be the least such. By 3(c), $n>0$. Assume also that there exists $m'>n$ such that $f_1(m')=f_2(m')$. Again, let $m$ be the least such. Then, $y_1=f_1(m-1)\neq y_2=f_2(m-1)$. Therefore, by 5(d), $R(y_1,\cdot)\cap R(y_2,\cdot)=\emptyset$ and by 5(e), $f_1(m)\in R(y_1,\cdot)$ and $f_2(m)\in R(y_2,\cdot)$ with $f_1(m)=f_2(m)$. Contradiction.
\end{proof}
\end{subclaim}

We will, now, define  $i:\M_1\cong\M_2$ by induction, so that, if $F(\M_1)=\{f_1,\ldots,f_n,\ldots\}$ and $F(\M_2)=\{g_1,\ldots,g_n,\ldots\}$, then $f_n$ (together with all the values $y=f_n(m)$) is included in the domain of $i$ at step $2n$, while at step $2n+1$ we make sure to include $g_n$ (and all the corresponding values). We do this in a way that the relations $R,V$ and the structures $M(y,\cdot)$ are preserved.

\textbf{Step $2n$:} Say that $i$ has been defined on \[X=\{f_1,\ldots, f_{n-1}, i^{-1}(g_1),\ldots, i^{-1}(g_{n-1})\},\] so far and that $f_n$ is not in $X$. Let $f\in X$ such that
\[m=s.p(f,f_n)\ge s.p(f',f_n),\mbox{ for all $f'\in X$.}\]
By inductive assumption, all the images (under $i$)  of the values $y_0=f(0)=f_n(0),\ldots,y_{m-1}=f(m-1)=f_n(m-1)$ have already  been defined. Let $z_j=i(y_j)$, for $j<m$. In $\M_2$ it holds by 5(e) that  $\bigwedge_{j<m-1} R(z_j,z_{j+1})$. By 5(b), there are infinitely many values in $R(z_{m-1},\cdot)$. Choose one of them, $z_m$,  that is different than all of $i(f_1)(m),\ldots, i(f_{n-1})(m)$ and $g_1(m),\ldots, g_{n-1}(m)$ and set $i(f_n(m))=z_m$. By 5(f), there is a function $g\in F(\M_2)$ that ``extends" $(z_0,\ldots,z_m)$. Let $i(f_n)=g$ and for $m'>m$, $i(f_n(m'))=g(m')$. Obviously $g$ is different than all of $i(f_1),\ldots, i(f_{n-1})$ and $g_1,\ldots, g_{n-1}$ and by 4(c) and 5, $R,V$ are preserved.

\textbf{Step $2n+1$:} Similarly.

Eventually, we will have included in the domain of $i$ the whole $F(\M_1)$ and in the range of $i$ all of $F(\M_2)$.  By subclaim \ref{3}, this also means that $V(n,\M_1), V(n,\M_2)$ are also included in the domain and the range of $i$ respectively.  As we mentioned $R$ and $V$ are preserved and $i$ becomes an isomorphism given that  we can extend $i$ on each of the structures $M(y,\cdot)$.

To this end, let $y\in \M_1$ and $z=i(y)$. By 6(c), $M(y,\M_1)\cup R(y,\M_1)\models \phi_{\lambda}$ and $M(z,\M_2)\cup R(z,\M_2)\models \phi_{\lambda}$. By completeness assumption on $\phi_{\lambda}$  and since both $\M_1,\M_2$ are countable, there is $j:M(y,\M_1)\cup R(y,\M_1)\cong M(z,\M_2)\cup R(z,\M_2)$. The problem is that $i,j$ may not agree on $R(y,\M_1)$. In either case, there exists a permuation $\pi$ of $R(y,\M_1)$ such that for all $y'\in R(y,\M_1)$, \[j(\pi(y'))=i(y').\]
By $P_{\lambda}$ being a homogeneous predicate, every such $\pi$ will induce an automorphism of $M(y,\M_1)\cup R(y,\M_1)$, call it $j_{\pi}$.  Then for every $y'\in R(y,\M_1)$  \[j(j_{\pi}(y'))=j(\pi(y'))=i(y').\]
So, we can extend $i$ on the whole of $M(y,\M_1)\cup R(y,\M_1)$ by $$i(y\rq{})=j(j_{\pi}(y\rq{}).$$
Since $j$ is onto $M(z,\M_2)\cup R(z,\M_2)$, we conclude that $i:\M_1\cong\M_2$.
\end{proof}
\end{claim}

The three previous claims complete the proof.
\end{proof}
\end{theorem}
Note: Obviously, if we could characterize $\omegaone$ by an \lomegaone- sentence, then we would also get characterizability of $\lambda^{\omegaone}$ etc. But this  is not possible. The inability  to characterize well-founded linear orderings also makes the obvious attempt to characterize $\ltok$ to fail.

We can actually do a bit better:

\begin{theorem}\label{ltoomegahom} If $\lambda\in\ch$, then $\ltoomega\in\homch$.
\begin{proof} First observe that  $\lambda$ is now characterizable and $\ltoomega$ is homogeneous characterizable. So, the assumption of the theorem is slightly weaker than the previous theorem and the conclusion is slightly stronger. 

If $\ltoomega\in\ch$, then $\ltoomega\in\homch$ by theorem \ref{kappatoomega} and it suffices to prove that $\ltoomega\in\ch$. We split the proof into two cases:
\begin{enumerate}
\item If $\ltoomega=\lambda$, then by theorem \ref{kappatoomega} again, $\lambda\in\homch$ and use the previous theorem.
\item If $\ltoomega\ge\lambda^+$, then, by Hausdorff formula
\[(\lambda^+)^{\omega}=\lambda^+\cdot\ltoomega=\ltoomega.\]

By theorem \ref{twohomcases}, one of $\lambda,\lambda^+$ is in $\homch$ and in either case we use the previous theorem.
\end{enumerate}
\end{proof}
\end{theorem}

\begin{corollary}\label{ltokhom} If $\kappa$ is an infinite cardinal and $\ltok\in\ch$, then $\ltok\in\homch$.
\end{corollary}

\begin{corollary}\label{chhomchcountableproducts} \ch and \homch\; are both closed under countable products.
\begin{proof} \[\prod_n {\alephs{\alpha_n}}=(\sup_n \alephs{\alpha_n})^{\alephs{0}}.\]
\end{proof}
\end{corollary}

\begin{theorem}\label{headpower} If $\alephalpha^{\alephs{\beta}}\in\ch$, then for all $\gamma<\omegaone$,
\[\alephs{\alpha+\gamma}^{\alephs{\beta}}\in\homch,\] i.e. if one power of $\alephalpha$ is in $\ch$, the same is true for the powers of a whole cluster of cardinals.
\begin{proof} First observe that by corollary \ref{ltokhom}, if $\alephs{\alpha+\gamma}^{\alephs{\beta}}$ is in $\ch$, then it is also in $\homch$. So, we do not worry about homogeneity. We proceed by induction on $\gamma$:

If $\gamma=\gamma_1+1$, a successor ordinal, then by the Hausdorff formula again
\[\alephs{\alpha+\gamma_1+1}^{\alephs{\beta}}=\alephs{\alpha+\gamma_1+1}\cdot\alephs{\alpha+\gamma_1}^{\alephs{\beta}}.\]
The second factor is in $\homch$ by inductive hypothesis and if  $\alephs{\alpha+\gamma_1+1}^{\alephs{\beta}}=\alephs{\alpha+\gamma_1}^{\alephs{\beta}}$, we can conclude that $\alephs{\alpha+\gamma_1+1}^{\alephs{\beta}}$ is also in $\homch$. If on the other hand $\alephs{\alpha+\gamma_1+1}^{\alephs{\beta}}>\alephs{\alpha+\gamma_1}^{\alephs{\beta}}$, then it must be the case that $\alephs{\alpha+\gamma_1+1}^{\alephs{\beta}}=\alephs{\alpha+\gamma_1+1}$ and  $\alephs{\alpha+\gamma_1}^{\alephs{\beta}}=\alephs{\alpha+\gamma_1}$. By the inductive hypothesis, $\alephs{\alpha+\gamma_1}^{\alephs{\beta}}$ is in $\ch$ and the same is true for $\alephs{\alpha+\gamma_1}$. We conclude that $\alephs{\alpha+\gamma_1+1}\in\ch$ by theorem \ref{successorthrm}.

If $\gamma=\lambda$, $\gamma\neq 0$, a countable limit ordinal,  then $cf(\alpha+\gamma)=cf(\gamma)=\omega$ and $\alephs{\alpha+\gamma}=\sup_n \alephs{\alpha+\gamma_n}$, for  an increasing sequence of $\gamma_n$'s. Then,
\[\alephs{\alpha+\gamma}^{\alephs{\beta}}=(\lim_n \alephs{\alpha+\gamma_n}^{\alephs{\beta}})^{cf(\alpha+\gamma)}=(\lim_n \alephs{\alpha+\gamma_n}^{\alephs{\beta}})^{\alephs{0}}=\prod_n (\alephs{\alpha+\gamma_n}^{\alephs{\beta}})\in\homch,\]
by the inductive assumption and corollary \ref{chhomchcountableproducts}.
\end{proof}
\end{theorem}

In particular, if we choose $\alpha=0$, then we conclude that for $\gamma<\omegaone$,

\[\alephs{\gamma}^{\alephs{\beta}}\in\homch\mbox{  iff  } 2^{\alephs{\beta}}\in\homch.\]
So, it is natural to ask

\textbf{Question 1:}  When does  it hold that $2^{\alephs{\beta}}\in\homch$?
In Part II we will consider ourselves with exactly this question. The importance of closure under the powerset  operation is also stressed by the next theorem.

\begin{theorem} \label{Cclass}If \C\; is the smallest set of characterizable cardinals that contains $\alephs{0}$
and is closed under successors, countable unions, countable products and powerset,
then it is also closed under powers.
\begin{proof}
First observe that
\begin{claim} If $\alephalpha^{\alephs{\beta}}\in C$, then for all $\gamma<\omegaone$,
$\alephs{\alpha+\gamma}^{\alephs{\beta}}\in C$.
\begin{proof} As in the proof of Theorem \ref{headpower} with the obvious modifications.
\end{proof}
\end{claim}

We will say that $\alephalpha$ is the head of a cluster in C (cf. definition \ref{headofcluster}), if there are no $\beta,\gamma$ such that
$\alephs{\beta}\in C$, $\gamma<\omegaone$ and $\alephalpha=\alephs{\beta+\gamma}$.

So, it suffices to prove that if $\alephalpha,\alephs{\beta}\in C$, then $\alephalpha^{\alephs{\beta}}\in C$, for all $\alephalpha$ heads of clusters.

If $\alephalpha\le\alephs{\beta}$, then $\alephalpha^{\alephs{\beta}}=2^{\alephs{\beta}}\in\C$. So,
assume that $\alephalpha>\alephs{\beta}\ge\alephs{0}$ and proceed by induction on $\alpha$.

Since $\alephalpha$ is the head of a cluster and by the way C was defined, $\alephalpha$ is either a countable union, or a countable product of smaller cardinals, or the powerset of a smaller cardinal.

\textbf{Case I} $\alephalpha=\sup_n \alephs{\alpha_n}$.

Then \[\alephalpha^{\alephs{\beta}}=(\lim_n \alephs{\alpha_n}^{\alephs{\beta}})^{\omega}.\]
For every $n$, $\alephs{\alpha_n}^{\alephs{\beta}}\in C$ by the inductive hypothesis and the result follows by closure under countable products. 

\textbf{Case II} $\alephalpha=\prod_{n} \alephs{\alpha_n}=(\sup_n \alephs{\alpha_n})^{\alephs{0}}$.

If $\alephalpha=\sup_n \alephs{\alpha_n}$, then we fall under Case I. So, assume that $\alephalpha>\sup_n \alephs{\alpha_n}$. Then
\[\alephalpha^{\alephs{\beta}}=(\sup_n \alephs{\alpha_n})^{\alephs{\beta}}\] and the result follows.

\textbf{Case III} $\alephalpha= 2^{\alephs{\gamma}}$, some $\alephs{\gamma}\in C$.

Then $\alephalpha^{\alephs{\beta}}=2^{\alephs{\gamma}\cdot\alephs{\beta}}$.
\end{proof}
\end{theorem}

\textbf{Question 2:} How much of the characterizable
cardinals this set $\C$ captures?

If we have GCH, then obviously
it captures everything. But are there any model of ZFC that give us examples of cardinals outside $\C$? We do not know the answer here. If not,
then we can tell all characterizable cardinals of a model quite
easily. In particular, the following conjecture of Shelah's follows easily and the answer is positive.
\begin{conjecture}[Shelah] If $\alephomegaone<2^{\alephs{0}}$, then every $\lomegaone$- sentence with  a model in power $\alephomegaone$ has a model in power  $2^{\alephs{0}}$.
\end{conjecture}
 If there are other cardinals outside $\C$, then it would be interesting to find such examples and even to  see how the whole picture of characterizable cardinals  looks like.

\section{Some counterexamples}
We provide some counterexamples to show that $\kappa\in\ch$
does not imply $cf(\kappa)\in\ch$, and $\kappa^+\in\ch$ doesn't
imply $\kappa\in\ch$.

\begin{theorem} \ch\; is not closed under \emph{predecessor} and \emph{cofinality}.
\begin{proof} We have to construct a counterexample. Shelah in \cite{ShelahBorelSquares} constructed a model where $2^{\alephs{0}}>\alephomegaone$ and where no cardinal $\kappa$ with $\alephomegaone\le\kappa< 2^{\alephs{0}}$ is characterizable. We say that in this case $\alephomegaone$ is the local Hanf number below $2^{\alephs{0}}$. This was done by adding \alephomegaone many Cohen reals in a ground model that satisfies GCH.

If we demand a little bit more here, we can take, for instance, that $2^{\alephs{0}}=\alephs{\omegaone+1}$. Then, $\alephs{\omegaone+1}$ is characterizable while \alephomegaone\; is not. This proves the first part of the theorem.

Now, if we let $2^{\alephs{0}}=\alephs{\omega_{\omegaone+1}}$, then we get \[cf(2^{\alephs{0}})= cf(\alephs{\omega_{\omegaone+1}})= \alephs{\omegaone+1},\] which, again, it is not a characterizable cardinal. This gives the second part of the theorem.
\end{proof}
\end{theorem}

\bibliographystyle{plain}
\bibliography{Bibliography}

\end{document}